\documentclass [a4paper,10pt]{article}
\usepackage [utf8] {inputenc}
\usepackage {amsmath}
\usepackage {amsthm}
\usepackage {amssymb}
\usepackage {amsfonts}
\newtheorem {lemma}{Lemma}
\newtheorem {theorem}{Theorem}
\title {A note on the distribution of $d_3(n)$ in arithmetic progressions}
\author {Tomos Parry}
\date {}
\begin {document}
\maketitle 
\begin {center}
\section {Introduction}
\end {center}
\begin {abstract}
Nguyen has shown that on averaging over $a=1,...,q$ the 3-fold divisor function has exponent of distribution 2/3, following \cite {banks}. We follow \cite {blomer} which leads to stronger bounds.
\end {abstract}
Let $d_k(n)$ be the $k$-fold divisor function
\[ d_k(n)=\sum _{a_1\cdot \cdot \cdot a_k=n}1.\]
What we would like is equidistribution over the primitive residue classes in the form 
\[ \sum _{n\leq x\atop {n\equiv a(q)}}d_k(n)=\text { main term }+\text { small error term }\]
where the main term is a simple function of size around $x/q$ and where the error term is a power-saving on the main term - crucial is the range of validity of $q$.
\\
\\ For $k=2$ we have a classical result of Selberg and of Hooley that says we have equidistribution for $q$ up to around $x^{2/3}$ and this seems to have resisted much improvement. Considering instead the average quantities
\[ \mathcal V_j'=\sideset {}{'}\sum _{a=1}^q\left |\sum _{n\leq x\atop {n\equiv a(q)}}d(n)-\text { main term }\right |^j\hspace {10mm}\mathcal V_j=\sum _{a=1}^q\left |\sum _{n\leq x\atop {n\equiv a(q)}}d(n)-\text { main term }\right |^j\]
where $\Sigma '$ means the sum is restricted to summing over residues coprime to the modulus.
Banks, Heath-Brown and Shparlinski \cite {banks} showed that averaging over $a$ is indeed enough to get a power-saving error term when $q$ goes pretty much all the way up to $x$. Specifically they showed
\begin {equation}\label {1}
\mathcal V_1'\ll x^\epsilon \left \{ \begin {array}{ll}
x^{3/4}&\text { if }q\leq x^{1/6}
\\ x^{2/3}q^{1/2}&\text { if }x^{1/6}<q\leq x^{1/3}
\\ x^{7/10}q^{2/5}&\text { if }x^{1/3}<q\leq x^{1/2}
\\ x^{4/5}q^{1/5}&\text { if }x^{1/2}<q\leq x
\end {array}\right .
\end {equation}
by using average bounds for incomplete Kloosterman sums. Blomer \cite {blomer} gave a simpler argument using the Voronoi summation formula that shows
\[ \mathcal V_2\ll x^{1+\epsilon }\hspace {10mm}(\text {and therefore }\mathcal V_1'\ll x^\epsilon \sqrt {xq})\]
and by seperating a diagonal term Lau and Zhao \cite {lauzhao} improved this argument to
\begin {eqnarray*}
\mathcal V_2&=&\left \{ \begin {array}{ll}\mathcal O\left (x^\epsilon \left ((xq)^{1/2}+q^{1/3}x^{2/3}\right )\right )&\text { if }q\leq x^{1/2}
\\ \text { main term }+\mathcal O\left (x^\epsilon \left (x^{5/6}q^{1/6}+x^{5/4}/q^{1/2}\right )\right )&\text { if }x^{1/2}<q\leq x.\end {array}\right .
\end {eqnarray*}
In this note we'll be interested in corresponding results for $k=3$, the interest in this case coming from the links to bounded gaps between primes - see, for example, Remark 1.4 of \cite {kowalski} or Remark 1 of \cite {nguyen1}.
\\
\\ In this case the work of Friedlander and Iwaniec in the 1980's and then improvements by Fouvry, Kowalski and Michel (FKM) in 2014 meant we can take prime $q$ up to $x^{1/2+1/46}$ for equidistribution, but as for $k=2$ we'd expect to be able to take larger $q$ if we average. Nguyen \cite {nguyen1} has shown that the argument of Banks, Heath-Brown, Shparlinski works for $k=3$ and obtained 
\begin {equation}\label {2}
\sideset {}{'}\sum _{a=1}^q\left |\sum _{n\leq x\atop {n\equiv a(q)}}d_3(n)-\text { main term }\right |\ll x^\epsilon \left \{ \begin {array}{ll}
x^{11/12}&\text { if }1\leq q\leq x^{1/4}
\\ x^{7/9}q^{1/2}&\text { if }x^{1/4}<q\leq x^{4/9}
\\ x&\text { if }x^{4/9}<q\leq x^{1/2}\hspace {5mm}\text {(this from the FKM result)}
\\ x^{5/6}q^{1/4}&\text { if }x^{1/2}<q\leq x^{2/3}
\end {array}\right .
\end {equation}
which in particular gives the nice statement that we get cancellation all the way up to $x^{2/3}$, parallel to cancellation up to $x$ in \cite {banks}. In the same way that \cite {blomer} improved \eqref {1}, in this note we'll be concerned with improving \eqref {2}. 
\\ \begin {theorem}\label {t1}
Let
\begin {eqnarray*}
E_{a/q}(s)&=&\sum _{n=1}^\infty \frac {d_3(n)}{n^s}e\left (\frac {na}{q}\right )\hspace {10mm}f_q(s)=Res_{s=1}\left \{ \frac {E_{a/q}(s)x^s}{s}\right \}
\hspace {10mm}\Delta (a/q)=\sum _{n\leq x}d_3(n)e\left (\frac {na}{q}\right )-f_q(s).
\end {eqnarray*}
Then
\begin {eqnarray*}
\sum _{a=1}^q\left |\Delta (a/q)\right |^2\ll x^{1+\epsilon }q^{3/2}.
\end {eqnarray*}
\end {theorem}
Let 
\[ \delta =(q,a)\hspace {10mm}\mathcal M_x(q,a)=\frac {x}{\phi (q/\delta )}Res_{s=1}\left \{ \left (\sum _{n=1\atop {(n,q)=\delta }}^\infty \frac {d_3(n)}{n^s}\right )\frac {x^{s-1}}{s}\right \} \hspace {10mm}E_x(q,a)=\sum _{n\leq x\atop {n\equiv a(q)}}d_3(n)-\mathcal M_x(q,a).\]
As
\begin {eqnarray*}
\sum _{a=1}^q|E_{x}(q,a)|^2&=&\frac {1}{q}\sum _{a=1}^q\left |\Delta (a/q)\right |^2
\end {eqnarray*}
the Cauchy-Schwarz Inequality with Theorem \ref {t1} gives
\begin {theorem}\label {t2}
Then, with $E_x(q,a)$ as above,
\begin {eqnarray*}
\sum _{a=1}^q|E_{x}(q,a)|\ll x^{1/2+\epsilon }q^{3/4}.
\end {eqnarray*}
\end {theorem}
Perhaps something stronger still could be obtained by exploiting cancellation within the Ramanujan sums at one stage in the proof (we have just taken absolute values of the dual sums). And one final comment: As remarked above the argument in \cite {blomer} and \cite {lauzhao} use the Voronoi summation formula, specifically the summation formula for $d(n)$. Ivi\' c has provided a Voronoi formula for $d_k(n)$, but this doesn't seem in the literature to have been applied too many times in work around equidistribution of the higher order divisor functions, despite their importance. With this result we are making it clear that this is possible and gets new results. There may be further applications too, and in a future paper we, for example, apply the formula to the study of the $L^1$-mean of the exponential sum of $d_3(n)$.

\begin {center}
\section {Lemmas}
\end {center}
First we state our main tool - Voronoi's summation formula for $d_3(n)$, due to Ivi\' c. Recall that for fixed $\sigma $ and $|t|\rightarrow \infty $ the Gamma function satisfies 
\[ |\Gamma (s)|\asymp |t|^{\sigma -1/2}e^{-\pi |t|/2}\]
and therefore
\begin {eqnarray*}
\left (\frac {\Gamma (s/2)}{\Gamma ((1-s)/2)}\right )^3\asymp \frac {1}{|t|^{3(1/2-\sigma )}}.
\end {eqnarray*}
In particular 
\[ \int _{(c)}\left (\frac {\Gamma (s/2)}{\Gamma ((1-s)/2)}\right )^3\frac {ds}{X^s}\]
converges absolutely for $0<c<1/6$ and any $X>0$.
\\ \begin {lemma}\label {voronoi3}
Let $E_{h/q}(s)$ be as in Theorem \ref {t1}, let $w:[0,\infty )\rightarrow \mathbb R$ be smooth and of compact support, and let
\begin {eqnarray*}
\tilde \Delta (h/q)&=&\sum _{n=1}^\infty d_3(n)e\left (\frac {nh}{q}\right )w(n)-Res_{s=1}\left (\int _0^\infty w(t)t^{s-1}dt\right )E_{h/q}(s)
\\U(X)&=&\frac {1}{2\pi i}\int _{(c)}\left (\frac {\Gamma (s/2)}{\Gamma ((1-s)/2)}\right )^3\frac {ds}{X^s}\hspace {10mm}\text {for }X>0\text { and }0<c<1/6
\\ N&=&\frac {\pi ^3n}{q^3}\hspace {20mm}\hat w_q(n)=\int _0^\infty w(t)U(Nt)dt
\\ R_{a,b,c}(h/q)&=&\sum _{x,y,z=1}^qe\left (\frac {ax+by+cz-hxyz}{q}\right )\hspace {10mm}A_{h/q}(n)=\sum _{abc=n}R_{a,b,c}(h/q).
\end {eqnarray*}
Then for $(h,q)=1$
\begin {eqnarray*}
\tilde \Delta (h/q)&=&\frac {\pi ^{3/2}}{q^3}\sum _{n=1}^\infty A_{h/q}(n)\hat w_q(n)+\text {\emph {3 similar terms}}.
\end {eqnarray*}
\end {lemma}
\begin {proof}
This is Theorem 2 of \cite {ivic} (the definitions of the various quantities in (2.1), (2.6), (2.7), (3.2) and line -4 of page 213, and the conditions on $f$ at the bottom of page 212). It is clear from the statement of that theorem that the similar terms have identical behaviour to the first.
\end {proof}
Let $c_q(n)$ denote Ramanujan's sum
\[ c_q(n):=\sideset {}{'}\sum _{a=1}^qe\left (\frac {an}{q}\right ).\]
We will use the following well-known properties of $c_q(n)$
\[ \text {for any $q,n,d\in \mathbb Z$}\hspace {10mm}\underbrace {\sum _{h|q,n}\mu (q/h)h}_{\ll q^\epsilon (q,n)}=c_{q}(n)=\frac {\phi (q)}{\phi (qd)}c_{qd}(nd)\hspace {10mm}\text { and }\hspace {10mm}(q,n)=1\implies c_q(n)=\mu (q)\]
without further mention. Denote by
\[ S_{n,m}(q)=\sideset {}{'}\sum _{a=1}^qe\left (\frac {na+m\overline a}{q}\right )\]
Kloosterman's sum. If $(n,q)=1$ then
\[ R_{a,b,c}(h/q)=qS_{1,\overline hn}(q)\]
and it's not too much trouble to then show
\[ \sideset {}{'}\sum _{h=1}^qA_{h/q}(n)\overline {A_{h/q}(m)}=q^3c_q(n-m)d_3(n)d_3(m)\]
but things get a bit involved once we don't have the coprimality condition. The next three lemmas deal with this.

\begin {lemma}\label {multiplicativity1}
Let $R_{a,b,c}(n)$ be as in Lemma \ref {voronoi3}. The sum
\begin {eqnarray*}
S_{\mathbf a,\mathbf b,\mathbf c}(q)=\sideset {}{'}\sum _{h=1}^qR_{a,b,c}(h/q)\overline {R_{a',b',c'}(h/q)}.
\end {eqnarray*}
is multiplicative.
\end {lemma}
\begin {proof}
The claim follows after establishing that, for $(h,q)=(h',q')=(q,q')=1$,
\begin {eqnarray*}
R_{a,b,c}\left (\frac {hq'+h'q}{qq'}\right )&=&R_{a,b,c}\left (\frac {hq'^3}{q}\right )R_{a,b,c}\left (\frac {h'q^3}{q'}\right ).
\end {eqnarray*}

\end {proof}

\begin {lemma}\label {R}
Suppose $q$ is a power of a prime $p$ and take $a,a',b,b'\in \mathbb Z$. Write
\begin {eqnarray*}
\mathcal B&=&(q,b,b')\hspace {10mm}Q=q/\mathcal B\hspace {10mm}\mathbf B=(b/\mathcal B,b'/\mathcal B)\hspace {10mm}\mathcal A=(q,a)\hspace {10mm}\mathcal A'=(q,a')\hspace {10mm}
\\ &&\hspace {30mm}\mathcal S=\sideset {}{'}\sum _{X,X'=1}^qe\left (\frac {aX-a'X'}{q}\right )c_q\left (bX-b'X'\right ).
\end {eqnarray*}
If $p|BB'$ then
\begin {eqnarray*}
\mathcal S&=&c_q(b,b')c_q(a)c_q(a').
\end {eqnarray*}
If $p\nmid BB'$ and $p^2|Q$ then
\begin {eqnarray*}
\mathcal S&=&qc_{q\mathcal B}(ab-a'b')\left \{ \begin {array}{ll}1&\text { if }\mathcal B=\mathcal A=\mathcal A'\\ 0&\text { if not}.\end {array}\right .
\end {eqnarray*}
If $p\nmid BB'$ and $Q=p$ then
\begin {eqnarray*}
\mathcal S&=&\underbrace {qc_{q\mathcal B}\left (ab'-a'b\right )}_{q/p|a,a'}-\mathcal Bc_q(a)c_q(a')
\end {eqnarray*}
\end {lemma}
In any case we have for $q$ a power of a prime $p$
\begin {eqnarray}\label {critical}
\sideset {}{'}\sum _{X,X'=1}^qe\left (\frac {aX-a'X'}{q}\right )c_q\left (bX-b'X'\right )&\ll &q\sum _{f|q(q,b,b'),ab-a'b'}f+q^\epsilon (q,b,b')\underbrace {(q,a)(q,a')}_{|q(q,a,a')}\notag 
\\ &\ll &q^{1+\epsilon }(q,b,b')\sum _{f|q,ab-a'b'}f.
\end {eqnarray}
\begin {proof}
We have
\begin {eqnarray}\label {dechra}
\mathcal S&=&\frac {\phi (q)}{\phi (Q)}\sideset {}{'}\sum _{X,X'=1}^qe\left (\frac {aX-a'X'}{q}\right )c_Q\left (BX-B'X'\right )\notag 
\\ &=&\frac {q}{Q}\sum _{h|Q}\mu (Q/h)h\sideset {}{'}\sum _{X,X'=1\atop {BX\equiv B'X'(h)}}^qe\left (\frac {aX-a'X'}{q}\right )\notag 
\\ &=:&\frac {q}{Q}\sum _{h|Q}\mu (Q/h)h\mathcal U(h)\notag 
\\ &=&q\left (\mathcal U(Q)-\frac {\mathcal U(Q/p)}{p}\right ).
\end {eqnarray}
If $p|BB'$ then $(p,BX-B'X')=1$ or $Q=1$ and the first equality above gives
\[ \mathcal S=\frac {\phi (q)\mu (Q)}{\phi (Q)}c_q(a)c_q(a')=c_q(\mathcal B)c_q(a)c_q(a')\]
which is the first claim, so let's suppose $p\nmid BB'$. For $1<h|Q$
\begin {eqnarray*}
\mathcal U(h)&=&\sideset {}{'}\sum _{X=1}^qe\left (\frac {aX}{q}\right )\sum _{X'=1\atop {(hX'+\overline {B'}BX,q)=1}}^{q/h}e\left (-\frac {a'}{q}\left (hX'+\overline {B'}BX\right )\right )
\\ &=&\frac {q}{h}\sideset {}{'}\sum _{X=1}^qe\left (\frac {X}{q}\left (a-a'\overline {B'}B\right )\right )\left \{ \begin {array}{ll}1&\text { if }q/h|a'\\ 0&\text { if not}\end {array}\right .
\\ &=&\frac {q}{h}c_q\left (aB'-a'B\right )\left \{ \begin {array}{ll}1&\text { if }q/h|a'\\ 0&\text { if not}\end {array}\right .
\end {eqnarray*}
so from \eqref {dechra} we have if $p^2|Q$
\begin {eqnarray*}
\mathcal S&=&qc_q\left (aB'-a'B\right )\left (\underbrace {\frac {1}{Q}}_{q/Q|a'}-\underbrace {\frac {1}{Q/p\cdot p}}_{qp/Q|a'}\right )
\\ &=&\frac {q\phi (q)}{Q\phi (q\mathcal B)}c_{q\mathcal B}\left (ab'-a'b\right )\left \{ \begin {array}{ll}1&\text { if }\mathcal B|a'\text { and }\mathcal Bp\nmid a'\\ 0&\text { if not}\end {array}\right .
\end {eqnarray*}
which gives the second claim (considering symmetry), and if $Q=p$ (meaning $\mathcal B=q/p$)
\begin {eqnarray*}
\mathcal S&=&q\left (\underbrace {\frac {q}{Q}c_q\left (aB'-a'B\right )}_{q/Q|a'}-\frac {c_q(a)c_q(a')}{p}\right )
=\underbrace {\frac {q^2\phi (q)}{Q\phi (q\mathcal B)}c_{q\mathcal B}(ab'-a'b)}_{q/Q|a'}-\mathcal Bc_q(a)c_q(a')
\end {eqnarray*}
which gives the third (again considering symmetry).

\end {proof}

\begin {lemma}\label {asum}
Let $R_{a,b,c}(h/q)$ be as in Lemma \ref {voronoi3} and write $n=abc$. Then
\begin {eqnarray*}
\sideset {}{'}\sum _{h=1}^qR_{a,b,c}(h/q)\overline {R_{a',b',c'}(h/q)}&\ll &q^{3+\epsilon }(q,n,n')\sum _{f|q,n-n'}f.
\end {eqnarray*}
\end {lemma}
\begin {proof}
By Lemma \ref {multiplicativity1} we can take $q$ a prime power. Since
\[ R_{a,b,c}(h/q)=q\sum _{\delta |q,c,b}\delta \sideset {}{'}\sum _{X=1}^{q/\delta }e\left (\frac {aX}{q/\delta }+\frac {bc/\delta ^2\overline h\overline X}{q/\delta }\right )\]
we have
\begin {eqnarray}\label {dimrhannu}
\text {for }q\nmid (b,c)\hspace {10mm}
R_{a,b,c}(h/q)&=&
q\sum _{\delta |q,c,b}
\sideset {}{'}\sum _{X=1}^{q}e\left (\frac {a\delta X+bc/\delta \overline h\overline X}{q}\right )\notag 
\\ \text {for }q|(b,c)\hspace {10mm}R_{a,b,c}(h/q)&=&q\sum _{d|q,a}d\phi (q/d)\leq q^{2+\epsilon }
\end {eqnarray}
so as long as we don't have ``$q|(b,c)$ and $q|(b',c')$" then the sum in question is 
\begin {eqnarray*}
\\ &=&q^2\sum _{\delta |q,c,b\atop {\delta '|q,c',b'}}\sideset {}{'}\sum _{X,X'=1}^{q}e\left (\frac {a\delta X-a'\delta 'X'}{q}\right )c_q\left (bc/\delta X'-b'c'/\delta 'X\right )
\\ &\ll &q^{3+\epsilon }(q,abc,a'b'c')\sum _{f|q,abc-a'b'c'}f\hspace {15mm}(\text {from the remark after Lemma \ref {R})}
\end {eqnarray*}
so we're ok in this case. If ``$q|(b,c)$ and $q|(b',c')$" (so $q|n,n'$) then \eqref {dimrhannu} says the LHS of the claim is $\ll q^{5+\epsilon }$ whilst the RHS is
\[ \geq q^{4+\epsilon }\sum _{f|q}f\]
so this case is also fine. Finally if we don't have ``$q|(b,c)$" but do have ``$q|(b',c')$"  (so $q|n'$) then \eqref {dimrhannu} say that the LHS of the claim is
\begin {eqnarray*}
&\leq &q^{3+\epsilon }\sum _{\delta |q,b,c}c_{q}(a\delta )c_q\left (bc/\delta \right )\ll q^{3+\epsilon }(q,n)
\end {eqnarray*}
whilst the RHS is
\begin {eqnarray*}
&\geq &q^{3+\epsilon }\sum _{f|q,n}f
\end {eqnarray*}
so we're good in all cases.
\end {proof}
Lemmas \ref {w1} and \ref {w2} are concerned with bounding the transform $\hat w_q(n)$ in the Voronoi summation formula.
\begin {lemma}\label {w1}
Let $1\leq Y\leq x$ and let $w:[0,\infty )\rightarrow \mathbb R$ be a smooth function satisfying
\begin {eqnarray*}
w(t)&=&\left \{ \begin {array}{ll}0&t\in [0,Y]\\ 1&t\in [2Y,x-Y]\\ 0&t\in [x,\infty )\end {array}\right .
\hspace {20mm}w^{(j)}(t)\ll \frac {1}{Y^j}\hspace {5mm}(j\geq 0)
\end {eqnarray*}
and let $U(X)$ be as in Lemma \ref {voronoi3}. 
If $Nx\gg 1$ then for any $j\in \mathbb N$
\begin {eqnarray*}
\int _0^\infty w(t)U\left (Nt\right )dt&\ll &\frac {Y}{(Nx)^{1/3}}\left (\frac {x^2}{NY^3}\right )^{j/3}.
\end {eqnarray*}
\end {lemma}
\begin {proof}
Let $S$ always denote any of the functions $\sin (6\cdot ^{1/3}),\cos (6\cdot ^{1/3})$ and let
\[ M(X)=\frac {S(X)}{X^{1/3}}\hspace {5mm}\text { and }\hspace {5mm}\mathcal M(N)=\int _0^\infty w(t)M\left (Nt\right )dt.\]
Then
\begin {eqnarray*}
N^{1/3}\mathcal M(N)
&=&\int _0^\infty \frac {w(t)S(Nt)}{t^{1/3}}dt\hspace {10mm}(\star )
\\ &=&\frac {-1}{2N^{1/3}}\int _0^\infty \frac {d}{dt}\left \{ w(t)t^{1/3}\right \} S(Nt)dt
\\ \text {so }\hspace {10mm}N^{1/3}|\mathcal M(N)|&\leq &\frac {1}{2N^{1/3}}\left (\underbrace {\left |\int _0^\infty \frac {w'(t)S(Nt)}{t^{1/3}}\cdot t^{2/3}dt\right |}_{\text {like $(\star )$ but multiplied by }t^{2/3}}+\underbrace {\left |\int _0^\infty \frac {w(t)S(Nt)}{t^{1/3}}\cdot \frac {dt}{t^{1/3}}\right |}_{\text {like }(\star )\text { but multiplied by $1/t^{1/3}$}}\right ).
\end {eqnarray*}
If we do this $j$ times we get
\begin {eqnarray*}
N^{1/3}|\mathcal M(N)|&\leq &\frac {1}{2^jN^{j/3}}\underbrace {\sum }_{2^j\text { terms, each with $A+B=j$}}\left |\int _0^\infty \frac {w^{(A)}(t)S(Nt)}{t^{1/3}}\cdot \frac {t^{2A/3}}{t^{B/3}}dt\right |
\\ &\ll &\frac {1}{2^jN^{j/3}}\underbrace {\sum \left (Y^{2/3-j/3}+Y^{1-A}x^{A-j/3-1/3}\right )}_{2^j\text { terms, each with $A+B=j$}}
\\ &\ll &\frac {Y}{x^{1/3}}\left (\frac {x^2}{NY^3}\right )^{j/3}
\end {eqnarray*}
which is the claim for $M(X)$ instead of $U(X)$. But from Lemma 3 of \cite {ivic} we have
\begin {eqnarray*}
U(X)&=&\text { main term }+\mathcal O\left (\frac {1}{X^{100}}\right )
\end {eqnarray*}
where the main term is $\ll 1$ sums of the form $M(X)$ (and of lower order). 
\end {proof}

\begin {lemma}\label {w2}
Let $\hat w_q(n)$ be as in Lemma \ref {voronoi3} and suppose $n\gg q^3/x$.
\\
\\ Then
\begin {eqnarray*}
\hat w_q(n)&\ll &\left \{ \begin {array}{ll}x^{1+\epsilon }&\\ x^{1/3}q^2/n^{2/3}.&\end {array}\right .
\end {eqnarray*}
Suppose 
\[ \frac {x^2q^3}{Y^3}\geq x^{1/100}=:x^{\delta }\]
and let $\epsilon >0$.  If $n>\left (x^2q^3/Y^3\right )^{1+\epsilon }$ then
\begin {eqnarray*}
\hat w_q(n)&\ll _\epsilon &\frac {1}{x^{100}n^{4/3}}.
\end {eqnarray*}

\end {lemma}
\begin {proof}
From Lemma \ref {w1} we have
\begin {eqnarray*}
\hat w_q(n)&\ll &\frac {Y}{(nx/q^3)^{1/3}}\left (\frac {x^2q^3}{nY^3}\right )^{j/3}
\end {eqnarray*}
which shows the first claim (the $\ll x$ bound following trivially from $U(X)\ll 1/X^\epsilon $ for $X>0$). We then continue, using the inequality assumed in the lemma,
\begin {eqnarray*}
\hat w_q(n)&\ll &Y\left (\frac {x^2q^3}{nY^3}\right )^{j/3}\ll \frac {x^7}{n^{4/3}}\left (\frac {x^2q^3}{nY^3}\right )^{(j-4)/3}
\leq \frac {x^7}{n^{4/3}}x^{-\delta \epsilon (j-4)/3} 
\end {eqnarray*}
so we get the second claim on taking $j$ sufficiently large.
\end {proof}
\begin {center}
\section {Proof of Theorem \ref {t1}}
\end {center}
Let $E_{a/q}(s),f_q(s),\Delta (a/q)$ be as in Theorem \ref {t1} and 
assume that all bounds can include a factor $x^\epsilon $ which we don't write in explicitly.
\\
\\ From (2.9) and (2.14) of \cite {ivic} $E_{h/q}(s)$ has Laurent expansion 
\[ \sum _{n\geq -3}a_{h/q}(n)(s-1)^n\]
where $a_{h/q}(n)\ll 1/q$, whilst
\begin {eqnarray*}
\frac {x^{s}}{s}&=&\int _0^\infty w(t)t^{s-1}dt+\mathcal O\left (Y|x^{s-1}|\right )
\end {eqnarray*}
and therefore
\begin {eqnarray*}
\Delta (a/q)&\ll &\overbrace {\sum _{n=1}^\infty d_k(n)w(n)e\left (\frac {na}{q}\right )-Res_{s=1}\left \{ E_{a/q}(s)\int _0^\infty w(t)t^{s-1}dt\right \} }^{=:\tilde \Delta (a/q)}
\\ &&+\sum _{n\text { in two intervals}\atop {\text {of length $Y$}}}d_k(n)e\left (\frac {na}{q}\right )\left (1-w(n)\right )+\mathcal O\left (\frac {Y}{q}\right )
\end {eqnarray*}
which gives
\begin {eqnarray*}
\sideset {}{'}\sum _{a=1}^q|\Delta (a/q)|^2&=&
\sideset {}{'}\sum _{a=1}^q|\tilde \Delta (a/q)|^2+q\sum _{n,m\text { in two intervals}\atop {\text {of length $Y$}\atop {n\equiv m(q)}}}d_k(n)d_k(m)+\frac {Y^2}{q}\ll \sideset {}{'}\sum _{a=1}^q|\tilde \Delta (a/q)|^2+Y^2
\end {eqnarray*}
and so
\begin {eqnarray}\label {balwnsabloda}
\sum _{a=1}^q|\Delta (a/q)|^2=\sum _{d|q}\sideset {}{'}\sum _{a=1}^d|\Delta (a/d)|^2\ll \sum _{d|q}\left (\sideset {}{'}\sum _{a=1}^d|\tilde \Delta (a/d)|^2+Y^2\right ).
\end {eqnarray} 
Let 
\[ M=\frac {1}{x}\left (\frac {xd}{Y}\right )^{3+\epsilon }.\]
From the first claim of Lemma \ref {w2}
\begin {eqnarray*}
\sum _{n\leq M\atop {n\equiv m(f)\atop {n\not =m}}}d_3(n)|\hat w_d(n)|&\ll &\sum _{n\ll d^3/x\atop {n\equiv m(f)\atop {n\not =m}}}x+\sum _{n\leq M\atop {n\equiv m(f)\atop {n\not =m}}}\frac {x^{1/3}d^2}{n^{2/3}}
\ll 
\frac {x^{1/3}d^2M^{1/3}}{f}
\end {eqnarray*}
so
\begin {eqnarray*}
\sum _{n,m\leq M\atop {n\equiv m(f)\atop {n\not =m}}}d_3(n)d_3(m)|\hat w_d(n)\hat w_d(m)|
\ll \frac {x^{2/3}d^4M^{2/3}}{f}
\ll \frac {x^2d^6}{Y^2f}
\end {eqnarray*}
as well as
\begin {eqnarray*}
\sum _{n,m\leq M\atop {n=m}}d_3(n)d_3(m)|\hat w_d(n)\hat w_d(m)|
&\ll &x^2\sum _{n\ll d^3/x}1+x\sum _{n\leq y}\frac {x^{1/3}d^2}{n^{2/3}}+\sum _{n>y}\frac {x^{2/3}d^4}{n^{4/3}}
\ll xd^3
\end {eqnarray*}
so, using the second claim of Lemma \ref {w2},
\begin {eqnarray*}
\sum _{n,m=1\atop {n\equiv m(f)}}^\infty d_3(n)d_3(m)|\hat w(n)\hat w(m)|&\ll &\frac {x^2d^6}{Y^2f}+xd^3+\text { tiny.}
\end {eqnarray*}
Therefore Lemmas \ref {voronoi3} and \ref {asum} give
\begin {eqnarray*}
\sideset {}{'}\sum _{h=1}^d|\tilde \Delta (h/d)|^2&=&\frac {\pi ^3}{d^6}\sum _{n,m=1}^\infty \hat w(n)\overline {\hat w(m)}\sideset {}{'}\sum _{h=1}^dA_{h/d}(n)\overline {A_{h/d}(m)}+\text { similar}
\\ &\ll &
\sum _{f|d}\frac {f}{d^3}\sum _{n,m=1\atop {n\equiv m(f)}}^\infty d_3(n)d_3(m)(d,n,n')|\hat w(n)\hat w(m)|
\\ &\ll &\left (\frac {xd}{Y}\right )^2d+xd
\end {eqnarray*}
and Theorem \ref {t1} follows from \eqref {balwnsabloda}, on choosing $Y=x^{1/2}q^{3/4}$.

\begin {center}

\end {center}

\end {document}